\documentclass[12pt]{article}

\usepackage[a4paper, mag=1000, includefoot,
            left=3cm, right=1cm, top=2cm, bottom=2cm, headsep=1cm, footskip=1cm]{geometry}
\usepackage[cp1251]{inputenc}
\usepackage[english,russian]{babel}
\usepackage{indentfirst}
\usepackage{footnpag}
\usepackage{amssymb} 
\usepackage{amsmath}
\usepackage{amsthm}
\usepackage{eucal}
\usepackage{dsfont}
\usepackage{graphicx}
\usepackage{psfrag}
\usepackage{caption}
\usepackage[numbers,square,comma,compress]{natbib}

\usepackage{ifpdf}
\ifpdf\usepackage{cmap}\fi 

\allowdisplaybreaks

\newtheorem{Def}{Определение}
\newtheorem{sta}{Утверждение}
\newtheorem*{teo*}{Теорема}
\newtheorem*{lem*}{Лемма}
\newtheorem*{con*}{Следствие}

\DeclareMathOperator{\re}{Re}
\DeclareMathOperator{\im}{Im}

\newcommand{\lt}{\left}
\newcommand{\rt}{\right}

\newcommand{\ts}{\times}

\renewcommand{\le}{\leqslant}
\renewcommand{\ge}{\geqslant}
\newcommand{\hm}[1]{#1 \nolinebreak \discretionary{}{\mbox{$#1$}}{}}

\newcommand{\al}{\alpha}
\newcommand{\bt}{\beta}
\newcommand{\gm}{{\gamma{}}}
\newcommand{\Gm}{{\Gamma{}}}

\newcommand{\vph}{\varphi}

\newcommand{\Om}{\Omega}
\newcommand{\la}{\lambda}
\newcommand{\La}{\Lambda}

\newcommand{\tl}{\tilde}

\newcommand{\ol}{\overline}

\newcommand{\RR}{\mathds{R}}
\newcommand{\NN}{\mathds{N}}

\newcommand{\CC}{\mathds{C}}

\begin{document}

\begin{otherlanguage}{english}

\title{On the Global Continuity of the Roots\\ of Families of Monic Polynomials (in Russian)\\ \selectlanguage{russian} О глобальной непрерывности корней\\ семейств приведённых полиномов}

\author{Evgeny E. Bukzhalev} 

\maketitle

\begin{abstract}

We raise a question on the existence of continuous roots of families of monic polynomials (by the root of a family of polynomials we mean a function of the coefficients of polynomials of a given family that maps each tuple of coefficients to a root of the polynomial with these coefficients). We prove that the family of monic second-degree polynomials with complex coefficients and the families of monic fourth-degree and fifth-degree polynomials with real coefficients have no continuous root. We also prove that the family of monic second-degree polynomials with real coefficients has continuous roots and we describe the set of all such roots.

\textbf{Keywords:} algebraic polynomial, polynomial function, continuous function, Jordan curve theorem, path-connected set, linear differential equation.

\end{abstract}

\end{otherlanguage}

\begin{abstract}

Поставлен вопрос о существовании непрерывных корней семейств приведённых многочленов (под корнем семейства многочленов понимается функция коэффициентов многочленов данного семейства, которая каждому набору коэффициентов ставит в соответствие один из корней многочлена с этими коэффициентами). Доказано, что семейство приведённых многочленов второй степени с комплексными коэффициентами и семейства приведённых многочленов четвёртой и пятой степеней с вещественными коэффициентами не имеют ни одного непрерывного корня. Также доказано, что семейство приведённых многочленов второй степени с вещественными коэффициентами имеет непрерывные корни и описано множество всех таких корней.

\textbf{Ключевые слова:} алгебраический полином, полиномиальная функция, непрерывная функция, теорема Жордана, линейно связное множество, линейное дифференциальное уравнение.

\end{abstract}

\section{Введение}

Будем рассматривать семейства $\mathfrak P^n_\Om$ приведённых многочленов степени $n \in \NN$ (то есть многочленов вида $z^n + a_{n-1}\, z^{n-1} + \dots + a_1\, z + a_0$) c коэффициентами из $\Om \in \mathfrak K := \{\RR, \CC\}$. При этом для каждого $(n,\Om) \in \NN \ts \mathfrak K$ под корнем семейства $\mathfrak P^n_\Om$ будем понимать такую функцию коэффициентов многочленов из $\mathfrak P^n_\Om$, которая каждому набору коэффициентов ставит в соответствие один из корней многочлена с этими коэффициентами. Основными результатами настоящей работы являются доказательства того, что семейства $\mathfrak P^2_\CC$, $\mathfrak P^4_\RR$ и $\mathfrak P^5_\RR$ не имеют ни одного непрерывного корня. Кроме того, мы также убедимся в том, что семейство $\mathfrak P^2_\RR$ имеет четыре непрерывных корня, и что эти корни образуют две пары полных множеств непрерывных корней этого семейства (см. определение \ref{def comp set}). Предложенные доказательства допускают обобщение на случай многочленов более высоких степеней --- с помощью аналогичных рассуждений можно доказать отсутствие непрерывных корней семейств $\mathfrak P^n_\CC$ при $n \ge 2$ и $\mathfrak P^n_\RR$ при $n \ge 4$.

Установленные результаты не противоречат утверждению о непрерывности корней семейства $\mathfrak P^n_\CC$ (при каждом натуральном $n$), доказанному в \cite{Sushkevich_41_book, Ostrowski_63_book, Voevodin_77_book, Harris_Martin_87_PAMS, Tyrtyshnikov_07_book1}, поскольку под непрерывностью во всех этих работах по существу подразумевается лишь непрерывность в любой наперёд заданной точке $(a_0, \dots, a_{n-1})$ пространства коэффициентов многочленов. Это утверждение может быть сформулировано следующим образом.
\begin{sta}\label{point cont}
  Для каждого натурального $n$ и каждой точки $M_0 \in \CC^n$ найдутся $n$ функций $z^1_{M_0}, \dots, z^n_{M_0}: \CC^n \to \CC$, непрерывных в точке $M_0$, и таких, что
  \[
    \mathfrak P^n_\CC(M) := z^n + a_{n-1}\, z^{n-1} + \dots + a_1\, z + a_0 = \big[ z - z^1_{M_0}(M) \big] \dots \big[ z - z^n_{M_0}(M) \big]
  \]
  при всех $M = (a_0, \dots, a_{n-1}) \in \CC^n$.
\end{sta}

Замечание. Каждое отображение $(z^1, \dots, z^n)$ непрерывно в своей точке $M_0$. То есть мы заранее фиксируем точку и под неё подбираем отображение, непрерывное в ней. При этом оно может быть разрывно в других точках. Мы можем сменить точку, но тогда, вообще говоря, придётся сменить и отображение. И не факт, что найдётся отображение, непрерывное сразу во всех точках. Таким образом, справедливость утверждения \ref{point cont} отнюдь не означает существования кортежа функций $(z^1, \dots, z^n)$ такого, что
\begin{gather}\label{tuple}
  \forall M \in \CC^n \quad \mathfrak P^n_\CC(M) = \big[ z - z^1(M) \big] \dots \big[ z - z^n(M) \big]
\end{gather}
(кортежи, удовлетворяющие условию \eqref{tuple}, далее будем называть полными наборами корней семейства $\mathfrak P^n_\CC$), каждая компонента $z^i$ которого непрерывна сразу во всём $\CC^n$. Всякая точка $M_0$ имеет своё множество $\mathfrak S_{M_0}$ полных наборов корней семейства $\mathfrak P^n_\CC$, непрерывных в этой точке (согласно утверждению \ref{point cont} каждое из $\mathfrak S_{M_0}$ заведомо не пусто), и пересечение семейства множеств $\{ \mathfrak S_{M_0} \}$ (по всем $M_0 \in \CC^n$) может не содержать ни одного кортежа.

Укажем на ещё одно свойство, связанное с непрерывностью корней семейств приведённых полиномов (доказательство см. в \cite{Tyrtyshnikov_07_book2}).
\begin{sta}\label{param pol}
  Рассмотрим отображение $\mathfrak p^n_{\mathds L}$, которое каждому $t \in \mathds L := [\al,\bt] \subset \RR$ ставит в соответствие многочлен $p^n_{\mathds L}(t) := x^{n} + a_{n-1}(t)\, x^{n-1} + \dots + a_1(t)\, x + a_0(t)$. Если $a_0, \dots, a_{n-1} \in C[\al,\bt]$, то существуют непрерывные функции $x_1, \dots, x_n \in C[\al,\bt]$ такие, что $p^n_{\mathds L}(t) = (x - x_1(t)) \dots (x - x_n(t))$.
\end{sta}

Отметим, что для справедливости данного утверждения существенна скалярность переменной $t$. В силу, например, утверждения \ref{pol 2C}, доказанного ниже, утверждение \ref{param pol} не может быть обобщено на случай векторной и/или комплексной переменной $t$.

Вопрос о непрерывности корней полиномов имеет связь с установлением равномерной экспоненциально-степенной оценки решения $w$ задачи Коши для линейного дифференциального уравнения с постоянными коэффициентами, рассматриваемыми как параметры для $w$:
\begin{gather}\label{de wn}
  \frac{d^nw}{{d\xi}^n} = a_{n-1}\, \frac{d^{n-1}w}{{d\xi}^{n-1}} + \cdots + a_0\, w, \quad \xi \in (0,+\infty);
\\\label{ic wn}
  w(0;M_n,N_n) = w^0,\ \ldots,\ \frac{d^{n-1}w}{{d\xi}^{n-1}}(0;M_n,N_n) = w^{n-1},
\end{gather}
где $M_n = (a_0, \dots, a_{n-1}) \in \CC^n$, $N_n = (w^0, \dots, w^{n-1}) \in \CC^n$. Справедливо утверждение (оно также содержит оценку производных $w$):

\begin{sta}

Пусть $\ol \La_n(M_n) := \max \{\re\la^1(M_n), \dots, \re\la^n(M_n)\}$, где $\la^i(M_n)$ --- корни характеристического многочлена $\la^n - a_{n-1}\, \la^{n-1} - \cdots - a_1\, \la - a_0$ уравнения \eqref{de wn} (подразумевается, что $(\la^1, \dots, \la^n)$ --- полный набор корней семейства характеристических многочленов по параметру $M_n = (a_0, \dots, a_{n-1}) \in \CC^n$), и пусть $\mathds D_a$ и $\mathds D_w$ --- замкнутые ограниченные подмножества пространства $\CC^n$. Тогда найдётся такое $\tl C_n \ge 0$, что
\begin{gather}\label{est w}
  \lt| \frac{d^iw}{{d\xi}^i}(\xi;M_n,N_n) \rt| \le \tl C_n\, (1 + \xi^{n-1})\, e^{\ol \La_n(M_n)\, \xi}
\end{gather}
при всех $(i,\xi,M_n,N_n) \in \{ 0, \dots, n-1 \} \ts [0, +\infty) \ts \mathds D_a \ts \mathds D_w$, где $w(\xi;M_n,N_n)$ --- решение задачи \eqref{de wn}--\eqref{ic wn}.

\end{sta}

Пусть известно (например, в силу критерия Гурвица --- см. \cite{Demidovich_67_book}), что $\forall (i,M_n) \in  \{ 0, \dots,$ $n-1 \} \ts \mathds D_a$ $\re\la^i(M_n) < 0$. Тогда $\forall M_n \in \mathds D_a$ $\ol \La_n(M_n) < 0$. Если все $\la^i$ непрерывны на $\mathds D_a$, то все $\re \la^i(M_n)$, а вслед за ними и $\ol \La_n$ (как максимум непрерывных функций), также непрерывны на $\mathds D_a$, а значит, в силу теоремы Вейерштрасса об экстремальных значениях непрерывной функции $\exists M^*_n \in \mathds D_a$ $\forall M_n \in \mathds D_a$ $\ol \La_n(M_n) \le \ol \La_n(M^*_n) < 0$. Таким образом, в случае непрерывных $\la^i$ сразу же имеем
\begin{gather*}
  \lt| \frac{d^iw}{{d\xi}^i}(\xi;M_n,N_n) \rt| \le \tl C_n\, (1 + \xi^{n-1})\, e^{- \varkappa\, \xi}
\end{gather*}
при всех $(i,\xi,M_n,N_n) \in \{ 0, \dots, n-1 \} \ts [0, +\infty) \ts \mathds D_a \ts \mathds D_w$, где $\varkappa := - \ol \La_n(M^*_n) > 0$.

Однако, как непосредственно следует из результатов настоящей статьи, семейство характеристических многочленов уравнения \eqref{de wn}, вообще говоря, не имеет полного набора корней, непрерывных на $D_a$ (более того, это семейство может не иметь даже одного такого корня), и поэтому справедливость равномерной экспоненциально-степенной оценки (с отрицательным коэффициентом в экспоненте) $\frac{d^i}{{d\xi}^i}\, w$ не настолько очевидна. Тем не менее, оказывается, что $\ol \La_n$ непрерывна во всём $\CC^n$ (а значит, и на любом наперёд заданном $D_a \subseteq \CC^n$) даже при разрывных $\la^i$ (для доказательства непрерывности $\ol \La_n$ в случае произвольных $\la^i$ проще всего воспользоваться утверждением \ref{point cont}, но, разумеется, это доказательство всё равно будет менее тривиальным, чем доказательство для непрерывных $\la^i$), и, тем самым, оценка \eqref{est w} в общем случае также остаётся в силе.

\section{Вспомогательный материал}

Пусть $n \in \NN_0 := \NN \cup \{0\}$, $\mathfrak P^n$ --- отображение, которое каждому $M = (a_0, \dots, a_{n-1}) \in \CC^n$ ставит в соответствие многочлен $\mathfrak P^n(M) := z^n + a_{n-1}\, z^{n-1} + \cdots + a_1\, z + a_0$. Поскольку разным точкам $M_1, M_2 \in \CC^n$ отвечают разные многочлены $\mathfrak P^n(M_1)$ и $\mathfrak P^n(M_2)$, то существует отображение ${(\mathfrak P^n)}^{-1}$, обратное к $\mathfrak P^n$, которое каждому многочлену $p \in \mathds P^n := \mathfrak P^n(\CC^n)$ ставит в соответствие такую точку $M \in \CC^n$, что $\mathfrak P^n(M) = p$ (здесь и ниже для всякой функции $\mathcal F: \mathcal X \to \mathcal Y$ и всякого множества $\mathcal M \subseteq \mathcal X$ под $\mathcal F(\mathcal M)$ подразумевается $\{y \in \mathcal Y \mid \exists x \in \mathcal M$ $\mathcal F(x) = y$).

Пусть $F^n: \CC^{n+1} \to \CC$, $(a_0, \dots, a_{n-1}, z) \mapsto z^n + a_{n-1}\, z^{n-1} + \cdots + a_1\, z + a_0$. Каждому $p \in \mathds P^n$ отвечает $F^n({(\mathfrak P^n)}^{-1}(p), \cdot): \CC \to \CC$, называемое полиномиальной функцией одного аргумента. Поскольку при любом $M \in \CC^n$ областью определения функции $F^n(M, \cdot)$ служит бесконечное множество $\CC$, то $\forall M_1, M_2 \in \CC^n\ M_1 \ne M_2 \Rightarrow F^n(M_1, \cdot) \ne F^n(M_2, \cdot)$. В силу последнего всякой функции $f \in \mathds F^n := F^n(\CC^n, \cdot)$ отвечает ровно одна точка $M \in \CC^n$ такая, что $F^n(M, \cdot) = f$, а значит, и ровно один многочлен $p \in \mathds P^n$ такой, что $F^n({(\mathfrak P^n)}^{-1}(p), \cdot) = f$.

Имея в виду установленное взаимно однозначное соответствие между множествами $\mathds P^n$ и $\mathds F^n$, вместо многочленов из $\mathds P^n$ будем рассматривать полиномиальные функции из $\mathds F^n$, которые в целях сокращения изложения будем называть просто многочленами. Вообще, далее под многочленом будет подразумеваться полиномиальная функция из $\mathds F^1 \cup \mathds F^2 \cup \dots$, под корнем многочлена $f$ --- нуль полиномиальной функции $f$, под коэффициентами многочлена $f$ --- координаты точки $M \in \CC^1 \cup \CC^2 \cup \dots$ такой, что $F^n(M, \cdot) = f$ (ввиду исключения из последующего рассмотрения обычных многочленов никакой двусмысленности при этом возникать не будет).

Пусть $\mathfrak F^n := \{(M,f) \in \CC^n \times \mathds F^n \mid f = F^n(M, \cdot)\}$ --- семейство многочленов $f \in \mathds F^n$ по параметру $M \in \CC^n$, или, что то же самое, $\mathfrak F^n: \CC^n \to \mathds F^n$, $M \mapsto F^n(M,\cdot)$. Заметим, что по определению образа множества $\mathds F^n = \mathfrak F^n(\CC^n)$. Для всякого $\Om \subseteq \CC$ обозначим $\mathfrak F^n_\Om$ сужение семейства $\mathfrak F^n$ на $n$-мерный куб $\Om^n$ (при этом $\mathfrak F^n_\CC = \mathfrak F^n$).

\begin{Def}
  Функция $r: \Om^n \to \CC$, называется корнем семейства многочленов $\mathfrak F^n_\Om$, если $r(M)$ является корнем многочлена $\mathfrak F^n_\Om(M)$ при каждом $M \in \Om^n$.
\end{Def}

Замечание 1. Во избежание недоразумений мы строго придерживаемся следующего принципа в обозначениях: $\mathcal F$ --- функция, $\mathcal F(x)$ --- её значение в точке $x$. Таким образом, если $r$ --- корень семейства многочленов, то есть отображение кортежа $M = (a_0,...,a_{n-1}) \hm\in \Om^n$ коэффициентов на один из корней многочлена c этими коэффициентами, то $r(M)$ --- один из корней многочлена с коэффициентами $M$.

Замечание 2. Поскольку под семейством многочленов мы понимаем отображение, которое каждому кортежу $M \in \Om^n$ коэффициентов ставит в соответствие многочлен c коэффициентами $M$, то в соответствии с принципом из замечания 1 мы обозначаем $\mathfrak F^n_\Om$ само семейство, а посредством $\mathfrak F^n_\Om(M)$ --- многочлен этого семейства (то есть значение отображения в точке $M$). Кроме того, мы различаем собственно полиномы (выражения) и полиномиальные функции. Понятно, что каждому полиному отвечает вполне определённая полиномиальная функция (и, с некоторыми оговорками, наоборот), и мы явно подчёркиваем это выше.

Далее нам понадобятся теорема Жордана и ряд связанных с ней определений (см., например, \cite{Spanier_71_book}).

\begin{Def}\label{path}
  Пусть $\al, \bt \in \RR$. Непрерывное отображение $\gm: [\al,\bt] \to \RR^2$ называется путём, соединяющим точки $\gm(\al)$ и $\gm(\bt)$.
\end{Def}

Замечание 3. Часто под путём подразумевают непрерывное отображение единичного отрезка (а не произвольного как в определении \ref{path}), однако привязка к единичному отрезку с технической точки зрения не всегда удобна, а замена в определении \ref{path} единичного отрезка произвольным существа определения \ref{path-con} (в котором используется понятие пути) не меняет.

\begin{Def}\label{path-con}
  Множество $\mathds M \subseteq \RR^2$ называется линейно связным в $\RR^2$, если для любых двух точек из $\mathds M$ существует путь $\gm$, соединяющий эти точки, образ которого принадлежит $\mathds M$.
\end{Def}

\begin{Def}
  Множество $\Gm \subseteq \RR^2$ называется кривой Жордана (или простой замкнутой кривой), если существует непрерывная инъекция $\gm: \Gm_0 \to \RR^2$ такая, что $\gm(\Gm_0) = \Gm$, где $\Gm_0$ --- единичная окружность на плоскости $\RR^2$.
\end{Def}

Замечание 4. Несложно доказать, что существование непрерывной инъекции $\gm: \Gm_0 \hm\to \RR^2$ такой, что $\gm(\Gm_0) = \Gm$, равносильно существованию $a,b \in \RR$ и вектор-функции $(\vph,\psi) \in C^2[a,b] := C[a,b] \times C[a,b]$ таких, что $a < b$, $\{(x,y) \in \RR^2 \mid \exists t \in [a,b]$ $(x = \vph(t)) \wedge (y \hm= \psi(t))\} = \Gm$, $(\vph(a),\psi(a)) = (\vph(b),\psi(b)) =: M_0$, $\forall M \in \Gm \backslash \{M_0\}$ $\exists! t \in (a,b)$ $(\vph(t),\psi(t)) = M$.

\begin{teo*}[Жордана]
  Пусть $\Gm$ есть простая замкнутая кривая на плоскости $\RR^2$. Тогда 1) множество $\mathds D := \RR^2 \backslash \Gm$ не является линейно связным в $\RR^2$, 2) существует единственная пара $\{\mathds D_1, \mathds D_2\}$ множеств, линейно связных в $\RR^2$, таких, что $\mathds D_1 \cup \mathds D_2 = \mathds D$ (при этом кривая $\Gm$ служит границей множеств $\mathds D_1$, $\mathds D_2$ и $\mathds D$).
\end{teo*}

Иными словами, любая простая замкнутая кривая на плоскости делит эту плоскость на две связные компоненты (согласно терминологии общей топологии $\mathds D_1$ и $\mathds D_2$ суть компоненты линейной связности множества $\mathds D$).

\section{Семейство $\mathfrak F^2_\RR$}

Пусть $\mathds S$ --- отображение, которое каждому $M \in \RR^2$ ставит в соответствие множество корней квадратичного полинома $\mathfrak F^2_\RR(M)$, $\mathfrak R := \{r: \RR^2 \to \CC \mid \forall M \in \RR^2$ $r(M) \in \mathds S(M)\}$ --- множество корней семейства $\mathfrak F^2_\RR$, $C^*(\Om)$ --- множество непрерывных функций из топологического пространства $\Om$ в пространство $\CC$, $\mathfrak R_C := \mathfrak R \cap C^*(\RR^2)$ --- множество непрерывных корней семейства $\mathfrak F^2_\RR$.

Пусть $\RR^2_\pm := \RR^2 \backslash \RR^2_0$, где $\RR^2_0 := \{(x,y) \in \RR^2 \mid 4\,x = y^2\}$ --- парабола на плоскости $\RR^2$, $\mathfrak X := 2^{\RR^2_\pm}$ (напомним, что для всякого множества $\mathcal M$ запись $2^{\mathcal M}$ означает булеан (множество всех подмножеств) множества $\mathcal M$), $S: \mathfrak X \ts \RR^2 \to \{-1,0,+1\}$,
\begin{gather*}
  (\mathds X, M) \mapsto
  \begin{cases}
    -1, & M \in \mathds X,
  \\
    \:\,0, & M \in \RR^2_0,
  \\
    +1, & M \in \RR^2_\pm \backslash \mathds X.
  \end{cases}
\end{gather*}

Пусть $\Xi$ --- отображение, которое каждому $\mathds X \in \mathfrak X$ ставит в соответствие функцию
\begin{gather}\label{roots}
  r: \RR^2 \to \CC,\ (a_0,a_1) \mapsto
  \begin{cases}
    {}-\dfrac{a_1}2 + S(\mathds X,(a_0,a_1))\, \dfrac{\sqrt{+D(a_0,a_1)}}2, & (a_0,a_1) \in \mathds \RR^2_+,
  \\[1ex]
    {}-\dfrac{a_1}2 + S(\mathds X,(a_0,a_1))\, \dfrac{\sqrt{-D(a_0,a_1)}}2\, i, & (a_0,a_1) \in \mathds \RR^2_-,
  \end{cases}
\end{gather}
где $D(a_0,a_1) := a_1^2 - 4\, a_0$ --- дискриминант полинома $\mathfrak F^2_\RR(a_0,a_1)$, $i$ --- мнимая единица, $\mathds \RR^2_+ :\hm= \{(x,y) \in \RR^2 \mid 4\,x < y^2\}$, $\mathds \RR^2_- := \{(x,y) \in \RR^2 \mid 4\,x > y^2\}$. Отображение $\Xi$ устанавливает взаимно однозначное соответствие между множествами $\mathfrak X$ и $\mathfrak R$ (это прямо следует из определений $\Xi$, $\mathfrak X$, $\mathfrak R$ и формул для корней квадратного трёхчлена).

Напомним, что $\forall M = (a_0, a_1) \in \RR^2$ $\mathfrak F^2_\RR(M)(z) := F^2(M,z) := z^2 + a_1\, z + a_0$.

\begin{Def}
  Двойка $R = (r_1,r_2) \in \mathfrak R^2$ называется полным набором корней семейства $\mathfrak F^2_\RR$, если $\forall (M,z) \in \RR^2 \ts \CC$ $F^2(M,z) = (z - r_1(M))(z - r_2(M))$. Если при этом $R \in \mathfrak R_C^2$, то $R$ называется полным набором непрерывных корней семейства $\mathfrak F^2_\RR$.
\end{Def}

\begin{Def}\label{def comp set}
  Пара $\mathcal R = \{\al,\bt\} \subseteq 2^{\mathfrak R}$ называется полным множеством корней семейства $\mathfrak F^2_\RR$, если $\forall (M,z) \in \RR^2 \ts \CC$ $F^2(M,z) = (z - \al(M))(z - \bt(M))$. Если при этом $\mathcal R \subseteq 2^{\mathfrak R_C}$, то $\mathcal R$ называется полным множеством непрерывных корней семейства $\mathfrak F^2_\RR$.
\end{Def}

Пусть $\mathfrak S$ --- множество полных наборов корней семейства $\mathfrak F^2_\RR$, $\mathfrak \Xi$ --- отображение, которое каждому $\mathds X \in \mathfrak X$ ставит в соответствие такую вектор-функцию $(r_1,r_2) \in {(\Xi(\mathfrak X))}^2$, что $r_1 \hm= \Xi(\mathds X)$, $r_2 = \Xi(\RR^2_\pm \backslash \mathds X)$. Отображение $\mathfrak \Xi$ устанавливает взаимно однозначное соответствие между множествами $\mathfrak X$ и $\mathfrak S$ (это прямо следует из определений $\mathfrak \Xi$, $\Xi$, $\mathfrak X$, $\mathfrak S$ и формул для корней квадратного трёхчлена). Заметим, что если $(r_1,r_2) \in \mathfrak S$, то и $(r_2,r_1) \in \mathfrak S$.

\begin{sta}\label{con root}
  Пусть $\mathfrak X_C := \{\varnothing, \RR^2_\pm, \RR^2_+, \RR^2_-\}$. Тогда $\mathfrak R_C = \Xi(\mathfrak X_C)$.
\end{sta}

\begin{proof}
Пусть $r \in \Xi(\mathfrak X_C)$. Убедимся, что $r \in \mathfrak R_C$. В силу определения $\Xi(\mathfrak X_C)$ найдётся $\mathds X_0 \in \mathfrak X_C$ такое, что $r = \Xi(\mathds X_0)$. Поскольку $\Xi$ --- биекция между множествами $\mathfrak X$ и $\mathfrak R$ и $\mathfrak X_C \subseteq \mathfrak X$, то $r \in \mathfrak R$. Остаётся доказать, что $r \in C^*(\RR^2)$.

Рассмотрим случай $\mathds X_0 = \RR^2_+$ (остальные три случая рассматриваются полностью аналогично). Так как $S(\RR^2_+,M) = -1$ при $M \in \RR^2_+$ и $S(\RR^2_+,M) = +1$ при $M \in \RR^2_-$, то сужения функции $r$ на области $\RR^2_+$ и $\RR^2_-$ принадлежат классу элементарных функций, и поэтому $r$ заведомо непрерывна в этих областях. Чтобы доказать непрерывность функции $r$ в произвольной точке $M_0 = (x_0,y_0)$ параболы $\RR^2_0$ (заметим, что $\RR^2 = \RR^2_+ \cup \RR^2_0 \cup \RR^2_-$), воспользуемся частичными пределами по множествам $\RR^2_+$, $\RR^2_0$ и $\RR^2_-$:
\begin{gather}\notag
  \lim_{\substack{(x,y) \to M_0\\ (x,y) \in \RR^2_+}} r(x,y) = \lim_{\substack{(x,y) \to M_0\\ (x,y) \in \RR^2_+}} \Bigg[ {}-\frac{y}2 + S(\RR^2_+,(x,y))\, \frac{\sqrt{D(x,y)}}2 \Bigg] =
\\\label{lim+}
  = \lim_{\substack{(x,y) \to M_0\\ (x,y) \in \RR^2_+}} \Bigg[ {}-\frac{y}2 - \frac{\sqrt{D(x,y)}}2 \Bigg] = {}-\frac{y_0}2 = r(x_0,y_0),
\end{gather}
и аналогично
\begin{gather}\label{lim0-}
  \lim_{\substack{(x,y) \to M_0\\ (x,y) \in \RR^2_0}} r(x,y) = r(x_0,y_0), \quad \lim_{\substack{(x,y) \to M_0\\ (x,y) \in \RR^2_-}} r(x,y) = r(x_0,y_0).
\end{gather}
Из \eqref{lim+} и \eqref{lim0-} следует, что $r(M) \to r(M_0)$ при $M \to M_0$, а значит, $r$ действительно непрерывна в точках $M_0$ из $\RR^2_0$.

Пусть теперь $r \in \mathfrak R_C$. Допустим, что $r \not\in \Xi(\mathfrak X_C)$. Поскольку $r \in \mathfrak R = \Xi(\mathfrak X)$ и $\mathfrak X_C \hm\subseteq \mathfrak X$, то найдётся $\mathds X_0 \in \mathfrak X \backslash \mathfrak X_C$ такое, что $r = \Xi(\mathds X_0)$. Заметим, что $\exists \chi \in \{\RR^2_+, \RR^2_-\}$ $\varnothing \hm\ne \chi \cap \mathds X_0 \ne \chi$ (иначе $\mathds X_0 \in \mathfrak X_C$). Предположим, что $\chi = \RR^2_+$ и положим $M_1 = (x_1,y_1) \hm\in \RR^2_+ \cap \mathds X_0$, $M_2(x_2,y_2) \in \RR^2_+ \backslash \mathds X_0$. Так как $r \in C^*(\RR^2)$, то из линейной связности множества $\RR^2_+$ (заметим также, что $r(\RR^2_+) = \RR$) вытекает линейная связность сужения
\[
  r_+ := \{(x,y,z) \in \RR^3 \mid ((x,y) \in \RR^2_+) \wedge (z = r(x,y))\}
\]
функции $r$ на $\RR^2_+$ (мы используем теоретико-множественный подход к понятию функции и отождествляем функцию $r_+: \RR^2_+ \to \RR$, $M \mapsto r(M)$ и её график). Из первой строки \eqref{roots} и определения $r_+$ видно, что $r_+ \subseteq r_1 \cup r_2$, где
\begin{align*}
  r_1 &:= \{(x,y,z) \in \RR^3 \mid ((x,y) \in \RR^2_+) \wedge (2\, z < -y)\},
\\
  r_2 &:= \{(x,y,z) \in \RR^3 \mid ((x,y) \in \RR^2_+) \wedge (2\, z > -y)\}.
\end{align*}
Для получения противоречия заметим, что (см. первую строку \eqref{roots})
\begin{gather*}
  r_+(M_1) := r(M_1) := {}-\frac{y_1}2 + S(\mathds X_0,M_1)\, \frac{\sqrt{D(M_1)}}2 = {}-\frac{y_1}2 - \frac{\sqrt{D(M_1)}}2 < {}-\frac{y_1}2,
\\
  r_+(M_2) := r(M_2) := {}-\frac{y_1}2 + S(\mathds X_0,M_2)\, \frac{\sqrt{D(M_2)}}2 = {}-\frac{y_1}2 + \frac{\sqrt{D(M_2)}}2 > {}-\frac{y_1}2,
\end{gather*}
а значит, $N_1 := (M_1,r_+(M_1)) \in r_+ \cap r_1$, $N_2 := (M_2,r_+(M_2)) \in r_+ \cap r_2$. Но $r_1$ и $r_2$ --- непересекающиеся области пространства $\RR^3$, поэтому существование точек $N_1$ и $N_2$ противоречит линейной связности $r_+$.

Предположение о том, что $\chi = \RR^2_-$, опровергается аналогично (при этом вместо функции $r$ удобнее рассмотреть $\im r$, воспользовавшись тем, что $\im r \in C(\RR^2)$). Таким образом, приходим к выводу, что $r \in \Xi(\mathfrak X_C)$.
\end{proof}

\begin{sta}
  Пусть $\mathfrak S_C$ --- множество полных наборов непрерывных корней семейства $\mathfrak F^2_\RR$. Тогда $\mathfrak S_C = \mathfrak \Xi(\mathfrak X_C)$.
\end{sta}

\begin{proof}
Пусть $R \in \mathfrak \Xi(\mathfrak X_C)$. Убедимся, что $R \in \mathfrak S_C = \mathfrak S \cap {(C^*(\RR^2))}^2$. В силу определения $\mathfrak \Xi(\mathfrak X_C)$ найдётся $\mathds X_0 \in \mathfrak X_C$ такое, что $R = \mathfrak \Xi(\mathds X_0) := (\Xi(\mathds X_0), \Xi(\RR^2_\pm \backslash \mathds X_0))$. Поскольку $\mathfrak \Xi$ --- биекция между множествами $\mathfrak X$ и $\mathfrak S$ и $\mathfrak X_C \subseteq \mathfrak X$, то $R \in \mathfrak S$. Для доказательства того, что $r \in {(C^*(\RR^2))}^2$, остаётся заметить, что $\forall \mathds X \in \mathfrak X_C$ $\RR^2_\pm \backslash \mathds X \in \mathfrak X_C$ и воспользоваться утверждением \ref{con root}.

Пусть $R = (r_1,r_2) \in \mathfrak S_C$. Убедимся, что $R \in \mathfrak \Xi(\mathfrak X_C)$. Поскольку $\mathfrak S_C = \mathfrak S \cap {(C^*(\RR^2))}^2$, а $\mathfrak S = \mathfrak \Xi(\mathfrak X)$, то найдётся $\mathds X_0 \in \mathfrak X$ такое, что $R = \mathfrak \Xi(\mathds X_0) := (\Xi(\mathds X_0), \Xi(\RR^2_\pm \backslash \mathds X_0))$. Так как при этом $r_1 = \Xi(\mathds X_0) \in \mathfrak R \cap C^*(\RR^2) =: \mathfrak R_C$, то для доказательства того, что $\mathds X_0 \in \mathfrak X_C$, остаётся ещё раз воспользоваться утверждением \ref{con root}.
\end{proof}

\begin{con*}
  Пусть $\mathcal S_C$ --- совокупность всех полных множеств непрерывных корней семейства $\mathfrak F^2_\RR$. Тогда $\mathcal S_C = \{\{\Xi(\varnothing), \Xi(\RR^2_\pm)\}, \{\Xi(\RR^2_+), \Xi(\RR^2_-)\}\}$.
\end{con*}

Замечание 1. С помощью формул Кардано для корней многочлена $x^3 + p\, x + q$ можно убедиться в существовании полного множества непрерывных корней семейства $\mathds P^3_\RR$, а также найти все отдельные непрерывные корни, все полные наборы и все полные множества непрерывных корней этого семейства.

Замечание 2. Анализ семейств $\mathfrak F^1_\RR$, $\mathfrak F^2_\RR$ и $\mathfrak F^3_\RR$ может привести к предположению о существовании полного множества корней (или хотя бы отдельного непрерывного корня) семейства $\mathfrak F^n_\RR$ при всяком натуральном $n$. Однако два последних раздела настоящей статьи опровергают данное предположение и показывают, что вещественность коэффициентов семейства приведённых многочленов степени $n$ не гарантирует существования даже одного непрерывного корня этого семейства.

\section{Семейство $\mathfrak F^2_\CC$}

Справедливо утверждение: при натуральном $n \ge 2$ семейство многочленов $\mathfrak F^n_\CC$ не имеет ни одного непрерывного корня (иными словами, при натуральном $n \ge 2$ не существует корня семейства приведённых многочленов степени $n$ с комплексными коэффициентами, являющегося непрерывной функцией коэффициентов многочленов этого семейства).

Мы ограничимся доказательством для случая $n=2$ (в случае произвольного натурального $n \ge 2$ доказательство более громоздко, но вполне аналогично).

\begin{sta}\label{pol 2C}
  Пусть $F: \CC^3 \to \CC,$ $(a_0, a_1, x) \mapsto x^2 + a_1\, x + a_0$. Тогда $\nexists r \in C^*(\CC^2)$ $\forall M \in \CC^2$ $F(M, r(M)) = 0$.
\end{sta}

$\triangle$ От обратного. Пусть $r \in C^*(\CC^2)$, $\forall M \in \CC^2$ $F(M, r(M)) = 0$. Обозначим $x_1(t) := e^{i\, \frac t2}$, $x_2(t) := -x_1(t)$, где $t \in \RR$. Раскрывая скобки, убеждаемся, что
\begin{gather*}
  (x - x_1(t))\, (x - x_2(t)) = x^2 - (x_1(t) + x_2(t))\, x + x_1(t)\, x_2(t) = x^2 - e^{i\,t} =
\\
  = x^2 + a_0(t) = F(a_0(t), 0, x) =: \tl F(t,x), \quad (t,x) \in \RR \ts \CC,
\end{gather*}
где $a_0(t) := -e^{i\,t} = {}- \cos t - i\, \sin t$.

Обозначим $\tl r(t) := r(a_0(t), 0)$, $\vph(t) := \re \tl r(t)$, $\psi(t) := \im \tl r(t)$, где $t \in \RR$. В силу теоремы о непрерывности сложной функции $\tl r \in C^*(\RR)$, $\vph, \psi \in C(\RR)$. Кроме того, очевидно, что $\tl r$, $\vph$ и $\psi$ суть $2\,\pi$-периодические функции переменной $t$. Зафиксируем два различных $t_1$ и $t_2$ из интервала $(0, 2\,\pi)$ и покажем, что $\tl r(t_1) \ne \tl r(t_2)$.

Поскольку согласно своему определению $\tl r(t)$ при каждом $t \in \RR$ является одним из корней многочлена $\tl F(t,\cdot)$, то существует функция $j: \RR \to \{1,2\}$ такая, что $\forall t \hm\in \RR$ $\tl r(t) = x_{j(t)}(t) =: x(j(t),t)$. Пусть $j(t_1) = j(t_2) =: j_0$. В силу формулы Эйлера для комплексной экспоненты из равенства $x(j_0,t_1) = x(j_0,t_2)$ следует, что $t_2 - t_1 = 4\, \pi\, k$ при некотором целом $k$. Но последнее противоречит тому, что $0 < |t_2 - t_1| < 2\, \pi$, поэтому $\tl r(t_1) = x(j_0,t_1) \ne x(j_0,t_2) = \tl r(t_2)$. Пусть теперь $j(t_1) \ne j(t_2)$. Вновь используя формулу Эйлера для комплексной экспоненты, при каждом $t \in (0, 2\,\pi)$ имеем $\im x_1(t) = \sin \frac t2 > 0$, $\im x_2(t) = -\sin \frac t2 < 0$. Значит, $\tl r(t_1) := x(j(t_1),t_1)$ и $\tl r(t_2) :\hm= x(j(t_2),t_2)$ опять заведомо не равны друг другу.

Так как $\tl r(0) = \tl r(2\,\pi)$, то из доказанного следует, что множество $\mathds L := \{ (x,y) \in \RR^2 \mid \exists t \hm\in [0, 2\,\pi]\ (x = \vph(t)) \wedge (y = \psi(t)) \}$ есть простая замкнутая кривая на плоскости $\RR^2$. Воспринимая комплексные числа как точки пространства $\RR^2$, введём в рассмотрение множество $\mathds M :\hm= \{ z \in \RR^2 \mid \exists k \in \{1,2\}$ $\exists t \in \RR$ $z = x_k(t) \}$. Из определений $x_k(t)$ видно, что $\mathds M \hm= \{(x,y) \in \RR^2 \mid x^2 + y^2 = 1\}$. Несложно доказать, что если $\varnothing \ne \mathds K \subseteq \mathds M$, то множество $(\RR^2 \backslash \mathds M) \cup \mathds K$ линейно связно в $\RR^2$.

Положим $\mathds K := \mathds M \backslash \mathds L$. При этом $\mathds L = \mathds M \backslash \mathds K$, так как из определений $\mathds M$ и $\mathds L$ вытекает, что $\mathds L \subseteq \mathds M$. Выше было установлено, что $\forall t \in (0, 2\,\pi)$ $\forall k \in \{1,2\}$ $\im x_k(t) \ne 0$, поэтому $\tl r(t) \notin \{ (-1,0), (+1,0) \} =: \mathds K_0$ при всех $t \in (0, 2\,\pi)$. Тогда, замечая, что
\[
  \tl r(0) = \tl r(2\,\pi) \in \bigcup_{k=1}^2 \{ x_k(0), x_k(2\,\pi) \} = \{ -1, +1 \} = \mathds K_0,
\]
приходим к выводу, что $\mathds L$ содержит ровно одну точку из $\mathds K_0$. При этом $\mathds K$ содержит оставшуюся точку из $\mathds K_0$, так как $\mathds K_0 \subseteq \mathds M$. Таким образом, в силу указанного выше признака множество $(\RR^2 \backslash \mathds M) \cup \mathds K = \RR^2 \backslash (\mathds M \backslash \mathds K) = \RR^2 \backslash \mathds L$ линейно связно в $\RR^2$. Но $L$~--- простая замкнутая кривая на плоскости $\RR^2$, а значит, последнее противоречит теореме Жордана, согласно которой множество $\mathds M \backslash \mathds L$ не является линейно связным на плоскости $\RR^2$. $\not\!\! \triangle$

Замечание 1. Ключевым (и, на мой взгляд, самым содержательным) моментом доказательства является установление непустоты множества $\mathds K$. Из неё следует, что множество $\mathds L$ не заполняет целиком окружность $\mathds M = \mathds L + \mathds K$ (являясь при этом её частью). Но если из окружности убрать множество, содержащее хотя бы одну точку $M$, то оставшееся множество --- в нашем случае множество $\mathds L$, не будет разделять плоскость $\mathds R^2$ на две не связанные компоненты (внутреннюю и внешнюю по отношению к окружности $\mathds M$ части плоскости $\mathds R^2$ можно будет соединить путём, <<проходящим>> через выброшенную из $\mathds M$ точку $M$). А должна делить в силу теоремы Жордана (так как по доказанному выше $\mathds L$~--- жорданова кривая). Полученное противоречие опровергает выдвинутое предположение о существовании непрерывного корня семейства $\mathfrak F^2_\CC$ и завершает доказательство теоремы.

Замечание 2. Используя римановы поверхности в качестве областей принадлежности коэффициентов семейства квадратичных полиномов, можно добиться существования непрерывного корня этого семейства. Здесь же мы, по сути, доказываем, что для обеспечения существования непрерывного корня семейства квадратичных полиномов c произвольными комплексными коэффициентами обойтись без многолистных римановых поверхностей невозможно: корень, определённый на $\CC^2$ --- прямом произведении двух однолистных поверхностей (обычных комплексных плоскостей), обязательно разрывен.

\section{Семейства $\mathfrak F^4_\RR$ и $\mathfrak F^5_\RR$}

Справедливо утверждение: при натуральном $n \ge 4$ семейство многочленов $\mathfrak F^n_{\RR}$ не имеет ни одного непрерывного корня (иными словами, при натуральном $n \ge 4$ не существует корня семейства приведённых многочленов степени $n$ с вещественными коэффициентами, являющегося непрерывной функцией коэффициентов многочленов этого семейства).

Мы ограничимся доказательством для случаев $n=4$ и $n=5$ (в случае произвольного натурального $n \ge 5$ доказательство более громоздко, но вполне аналогично доказательству для $n=5$).

\begin{sta}
  Пусть $F: \CC^5 \to \CC,$ $(a_0, a_1, a_2, a_3, x) \mapsto x^4 + a_3\, x^3 + a_2\, x^2 + a_1\, x + a_0$. Тогда $\nexists r \in C^*(\RR^4)$ $\forall M \in \RR^4$ $F(M, r(M)) = 0$.
\end{sta}

$\triangle$ От обратного. Пусть $r \in C^*(\RR^4)$, $\forall M \in \RR^4$ $F(M, r(M)) = 0$. Обозначим $x_1(t) :\hm= i + e^{i\, (\frac t2 - \frac\pi2)}$, $x_2(t) := i - e^{i\, (\frac t2 - \frac\pi2)}$, $x_3(t) := \bar x_1(t)$, $x_4(t) := \bar x_2(t)$, где $t \in \RR$, а черта означает комплексное сопряжение. Раскрывая скобки, убеждаемся, что
\begin{gather*}
  (x - x_1(t))\, (x - x_2(t)) = x^2 - (x_1(t) + x_2(t))\, x + x_1(t)\, x_2(t) = x^2 - 2\, i\, x - 1 + e^{i\,t},
\\
  (x - x_1(t))\, (x - x_2(t))\, (x - x_3(t))\, (x - x_4(t)) = (x^2 - 1 + \cos t - i\, (2\, x - \sin t)) \ts{}
\\
  {}\ts (x^2 - 1 + \cos t + i\, (2\, x - \sin t)) = {(x^2 - 1 + \cos t)}^2 + {(2\, x - \sin t)}^2 =
\\
  = x^4 + a_2(t)\, x^2 + a_1(t)\, x + a_0(t) = F(a_0(t), a_1(t), a_2(t), 0, x) =: \tl F(t,x), \quad (t,x) \in \RR \ts \CC,
\end{gather*}
где $a_0(t) := 2\, (1 - \cos t)$, $a_1(t) := - 4\, \sin t$, $a_2(t) := 2\, (1 + \cos t)$.

Обозначим $\tl r(t) := r(a_0(t), a_1(t), a_2(t), 0)$, $\vph(t) := \re \tl r(t)$, $\psi(t) := \im \tl r(t)$, где $t \in \RR$. В силу теоремы о непрерывности сложной функции $\tl r \in C^*(\RR)$, $\vph, \psi \in C(\RR)$. Кроме того, очевидно, что $\tl r$, $\vph$ и $\psi$ суть $2\,\pi$-периодические функции переменной $t$. Зафиксируем два различных $t_1$ и $t_2$ из интервала $(0, 2\,\pi)$ и покажем, что $\tl r(t_1) \ne \tl r(t_2)$.

Поскольку согласно своему определению $\tl r(t)$ при каждом $t \in \RR$ является одним из корней многочлена $\tl F(t,\cdot)$, то существует функция $j: \RR \to \{1,2,3,4\}$ такая, что $\forall t \hm\in \RR$ $\tl r(t) = x_{j(t)}(t) =: x(j(t),t)$. Пусть $j(t_1) = j(t_2) =: j_0$. В силу формулы Эйлера для комплексной экспоненты из равенства $x(j_0,t_1) = x(j_0,t_2)$ следует, что $t_2 - t_1 = 4\, \pi\, k$ при некотором целом $k$. Но последнее противоречит тому, что $0 < |t_2 - t_1| < 2\, \pi$, поэтому $\tl r(t_1) = x(j_0,t_1) \ne x(j_0,t_2) = \tl r(t_2)$. Пусть теперь $j(t_1) \ne j(t_2)$. Вновь используя формулу Эйлера для комплексной экспоненты, при каждом $t \in (0, 2\,\pi)$ имеем $\re x_1(t) = \cos \big( \frac t2 - \frac\pi2 \big) \hm= \sin \frac t2 > 0$, $\im x_1(t) = 1 + \sin \big( \frac t2 - \frac\pi2 \big) = 1 - \cos \frac t2 > 0$, и аналогично $\re x_2(t) < 0$, $\im x_2(t) > 0$, $\re x_3(t) > 0$, $\im x_3(t) < 0$, $\re x_4(t) < 0$, $\im x_4(t) < 0$. Но тогда либо $\re x(j(t_1),t_1) \hm\ne \re x(j(t_2),t_2)$, либо $\im x(j(t_1),t_1) \ne \im x(j(t_2),t_2)$, а значит, $\tl r(t_1) := x(j(t_1),t_1)$ и $\tl r(t_2) :\hm= x(j(t_2),t_2)$ опять заведомо не равны друг другу.

Так как $\tl r(0) = \tl r(2\,\pi)$, то из доказанного следует, что множество $\mathds L := \{ (x,y) \in \RR^2 \mid \exists t \hm\in [0, 2\,\pi]\ (x = \vph(t)) \wedge (y = \psi(t)) \}$ есть простая замкнутая кривая на плоскости $\RR^2$. Воспринимая комплексные числа как точки пространства $\RR^2$, введём в рассмотрение множество $\mathds M :\hm= \{ z \in \RR^2 \mid \exists k \in \{1,2,3,4\}$ $\exists t \in \RR$ $z = x_k(t) \}$. Из определений $x_k(t)$ видно, что $\mathds M \hm= \mathds M_1 \cup \mathds M_2$, где $\mathds M_1 := \{(x,y) \in \RR^2 \mid x^2 + (y-1)^2 = 1\}$, $\mathds M_2 := \{(x,y) \in \RR^2 \mid x^2 + (y+1)^2 = 1\}$. Несложно доказать, что если $\mathds K \subseteq \mathds M$ и $(\exists z_1 \in \mathds K \cap \mathds M_1) \wedge (\exists z_2 \in \mathds K \cap \mathds M_2)$, то множество $(\RR^2 \backslash \mathds M) \cup \mathds K$ линейно связно в $\RR^2$.

Положим $\mathds K := \mathds M \backslash \mathds L$. При этом $\mathds L = \mathds M \backslash \mathds K$, так как из определений $\mathds M$ и $\mathds L$ вытекает, что $\mathds L \subseteq \mathds M$. Выше было установлено, что $\forall t \in (0, 2\,\pi)$ $\forall k \in \{1,2,3,4\}$ $\re x_k(t) \cdot \im x_k(t) \ne 0$, поэтому $\tl r(t) \notin \{ (0,-2), (0,0), (0,+2) \} =: \mathds K_0$ при всех $t \in (0, 2\,\pi)$. Тогда, замечая, что
\[
  \tl r(0) = \tl r(2\,\pi) \in \bigcup_{k=1}^4 \{ x_k(0), x_k(2\,\pi) \} = \{ -2\, i, 0, +2\, i \} = \mathds K_0,
\]
приходим к выводу, что $\mathds L$ содержит ровно одну точку из $\mathds K_0$. При этом $\mathds K$ содержит две оставшиеся точки из $\mathds K_0$, так как $\mathds K_0 \subseteq \mathds M$. Но поскольку каждое из множеств $\mathds M_1$ и $\mathds M_2$ содержит по крайней мере по одной точке из любого двухточечного подмножества $\mathds K_0$, то в силу указанного выше признака множество $(\RR^2 \backslash \mathds M) \cup \mathds K = \RR^2 \backslash (\mathds M \backslash \mathds K) = \RR^2 \backslash \mathds L$ линейно связно в $\RR^2$. Однако $L$ --- простая замкнутая кривая на плоскости $\RR^2$, а значит, последнее противоречит теореме Жордана, согласно которой множество $\mathds M \backslash \mathds L$ не является линейно связным на плоскости $\RR^2$. $\not\!\! \triangle$

Замечание. Предыдущие два утверждения могут быть доказаны без применения теоремы Жордана. Альтернативный подход будет использован при доказательстве утверждения для $n=5$. Но схема доказательства утверждения в случае $n=5$, не предполагающая обязательного применения теоремы Жордана, может быть легко перенесена (с соответствующим упрощением) на случай $n=4$ для многочленов с коэффициентами из $\RR$ и случай $n=2$ для многочленов с коэффициентами из $\CC$. Использование альтернативного похода для $n=5$, разумеется, не является необходимым --- можно было ограничиться небольшой модификацией схемы (добавлением корня, <<скользящего>> вдоль действительной оси), использовавшейся при доказательстве предыдущих утверждений.

\begin{sta}
  Пусть $F: \CC^6 \to \CC,$ $(a_0, a_1, a_2, a_3, a_4, x) \mapsto x^5 + a_4\, x^4 + a_3\, x^3 + a_2\, x^2 + a_1\, x \hm+ a_0$. Тогда $\nexists r \in C^*(\RR^5)$ $\forall M \in \RR^5$ $F(M, r(M)) = 0$.
\end{sta}

$\triangle$ От обратного. Пусть $r \in C^*(\RR^5)$, $\forall M \in \RR^5$ $F(M, r(M)) = 0$. Положим
\begin{gather*}
  x_1(t) :=
  \begin{cases}
    i + e^{i\, (\frac t2 - \frac\pi2)}, & t \in [0, 2\,\pi],
  \\
    2\, i, & t \in [2\,\pi, 4\,\pi],
  \\
    i - e^{i\, (\frac t2 - \frac\pi2)}, & t \in [4\,\pi, 6\,\pi];
  \end{cases}
  \quad x_2(t) :=
  \begin{cases}
    i - e^{i\, (\frac t2 - \frac\pi2)}, & t \in [0, 2\,\pi],
  \\
    0, & t \in [2\,\pi, 4\,\pi],
  \\
    i + e^{i\, (\frac t2 - \frac\pi2)}, & t \in [4\,\pi, 6\,\pi];
  \end{cases}
\\
  x_3(t) := \bar x_1(t),\ t \in [0, 6\,\pi];\ x_4(t) :=
  \begin{cases}
    \bar x_2(t), & t \in [0, 2\,\pi],
  \\
    t - 2\,\pi, & t \in [2\,\pi, 4\,\pi],
  \\
    +2\,\pi, & t \in [4\,\pi, 6\,\pi];
  \end{cases}
  \ x_5(t) :=
  \begin{cases}
    -2\,\pi, & t \in [0, 2\,\pi],
  \\
    t - 4\,\pi, & t \in [2\,\pi, 4\,\pi],
  \\
    \bar x_2(t), & t \in [4\,\pi, 6\,\pi];
  \end{cases}
\\
  \tl F(t,x) := (x - x_1(t))\, (x - x_2(t))\, (x - x_3(t))\, (x - x_4(t))\, (x - x_5(t)) =
\\
  = x^5 + a_4(t)\, x^4 + a_3(t)\, x^3 + a_2(t)\, x^2 + a_1(t)\, x + a_0(t), \quad \tl r(t) := r(a_0(t), a_1(t), a_2(t), a_3(t), a_4(t)),
\end{gather*}
где $a_k: [0, 6\, \pi] \to \RR$, $t \mapsto \tl a_k(x_1(t),x_2(t),x_3(t),x_4(t),x_5(t))$, $\tl a_k$ --- полиномиальные (а значит и непрерывные) функции пяти аргументов. Вещественнозначность функций $a_k$ является следствием того, что при каждом $t \in [0, 6\, \pi]$ множество $\{ z \in \CC \backslash \RR \mid \exists k \in \NN_5$ $z = x_i(t)\}$ (здесь и далее $\NN_5 := \{1,2,3,4,5\}$) разбивается на комплексно сопряжённые пары. Используя теорему о непрерывности сложной функции, последовательно заключаем, что все $a_k$ и $\tl r$ непрерывны на отрезке $[0, 6\,\pi]$.

Получим точные выражения для $a_k(t)$ при $t \in T := [0, 2\,\pi] \cup [2\,\pi, 6\,\pi]$. Раскрывая скобки и считая, что $t \in T$, убеждаемся, что
\begin{gather*}
  (x - x_1(t))\, (x - x_2(t)) = x^2 - (x_1(t) + x_2(t))\, x + x_1(t)\, x_2(t) = x^2 - 2\, i\, x - 1 + e^{i\,t},
\\
  (x - x_1(t))\, (x - x_2(t))\, (x - \bar x_1(t))\, (x - \bar x_2(t)) = (x^2 - 1 + \cos t - i\, (2\, x - \sin t)) \ts{}
\\
  {}\ts (x^2 - 1 + \cos t + i\, (2\, x - \sin t)) = {(x^2 - 1 + \cos t)}^2 + {(2\, x - \sin t)}^2 =
\\
  = x^4 + b_2(t)\, x^2 + b_1(t)\, x + b_0(t) =: \tl F_4(t,x),
\end{gather*}
где $b_0(t) := 2\, (1 - \cos t)$, $b_1(t) := - 4\, \sin t$, $b_2(t) := 2\, (1 + \cos t)$. Тогда
\begin{gather*}
  \tl F(t,x) =
  \begin{cases}
    \tl F_4(t,x)\, (x - x_5(t)) = \tl F_4(t,x)\, (x + 2\,\pi), & t \in [0, 2\,\pi],
  \\
    \tl F_4(t,x)\, (x - x_4(t)) = \tl F_4(t,x)\, (x - 2\,\pi), & t \in [4\,\pi, 6\,\pi].
  \end{cases}
\end{gather*}
Отсюда, вновь раскрывая скобки и приводя подобные при $x$, имеем: $a_4(t) = s(t)$, $a_3(t) \hm= b_2(t)$, $a_2(t) = b_1(t) + s(t)\, b_2(t)$, $a_1(t) = b_0(t) + s(t)\, b_1(t)$, $a_0(t) = s(t)\, b_0(t)$, где $s(t) := +2\,\pi$ при $t \in [0, 2\,\pi]$ и $s(t) := -2\,\pi$ при $t \in [4\,\pi, 6\,\pi]$.

Из точных выражений для $a_k(t)$ видно, что $a_k(0) = a_k(2\,\pi)$, $a_k(4\,\pi) = a_k(6\,\pi)$ при всех $k \in \{0,1,2,3,4\}$.

Будем интерпретировать комплексные числа как точки пространства $\RR^2$. Поскольку согласно своему определению $\tl r(t)$ при каждом $t \in \RR$ является одним из корней $x_j(t) \hm=: x(j,t)$ многочлена $\tl F(t,\cdot)$, то\footnote{В случае многочлена произвольной степени условие принадлежности части окружности можно задать двумя строгими неравенствами (знаками отклонений), получаемыми из нормальных уравнения прямых.}
\begin{gather*}
  \tl r((0, 2\,\pi)) \subseteq X_1 := x(\NN_5 \ts (0, 2\,\pi)) = \bigcup_{k=1}^5 x_k(0, 2\,\pi) = \{(x,y) \in \Om_1(0,+1) \mid x > 0\} \cup{}
\\
  {}\cup \{(x,y) \in \Om_1(0,+1) \mid x < 0\} \cup \{(x,y) \in \Om_1(0,-1) \mid x > 0\} \cup{}
\\
  {}\cup \{(x,y) \in \Om_1(0,-1) \mid x < 0\} \cup \{(-2\,\pi, 0)\},
\end{gather*}
где $\Om_R(A)$ --- окружность радиуса $R$ c центром в т. $A$ на плоскости $\RR^2$. Несложно видеть\footnote{Для строго обоснования, на мой взгляд, проще всего воспользоваться всё той же теоремой Жордана.}, что $x_k((0, 2\,\pi))$ суть компоненты линейной связности множества $X_1$ (это значит, что каждое из $x_k((0, 2\,\pi))$ линейно связно в $\RR^2$, но если к какому-нибудь из $x_k((0, 2\,\pi))$ добавить хотя бы одну точку из $X_1 \backslash x_k((0, 2\,\pi))$, то получившееся множество уже не будет линейно связным в $\RR^2$). Но при непрерывном отображении образ линейно связного множества (в данном случае интервала $(0, 2\,\pi)$) линейно связен, поэтому $\tl r((0, 2\,\pi)) \subseteq x_k((0, 2\,\pi))$ при некотором $k \in \NN_5$.

Пусть $\tl r((0, 2\,\pi)) \subseteq x_1((0, 2\,\pi))$. Так как при этом $x_1((0, 2\,\pi))$ не пересекается с другими $x_k((0, 2\,\pi))$, то $\forall t \in (0, 2\,\pi)$ $\tl r(t) = x_1(t)$. Тогда, ещё раз вспоминая о том, что $\tl r$ непрерывна на отрезке $[0, 6\pi]$ (а значит, непрерывна во всех точках этого отреза), имеем:
\[
  \tl r(0) = \lim_{t \to +0} \tl r(t) = \lim_{t \to +0} x_1(t) = 0, \quad \tl r(2\,\pi) = \lim_{t \to 2\pi-0} \tl r(t) = \lim_{t \to 2\pi-0} x_1(t) = 2\, i.
\]
Но $\tl r(0) := r(a_0(0), \dots, a_4(0)) = r(a_0(2\,\pi), \dots, a_4(2\,\pi)) =: \tl r(2\,\pi)$. Следовательно, $\tl r((0, 2\,\pi)) \hm\nsubseteq x_1((0, 2\,\pi))$. Предположения о том, что $\tl r((0, 2\,\pi)) \subseteq x_k((0, 2\,\pi))$ при $k \in \{2,3,4\}$, опровергаются полностью аналогично.

Остаётся только одна возможность: $\tl r((0, 2\,\pi)) \subseteq x_5((0, 2\,\pi))$. Так как $x_5((0, 2\,\pi))$ также не пересекается с другими $x_k((0, 2\,\pi))$, то $\forall t \in (0, 2\,\pi)$ $\tl r(t) = x_5(t)$. Отсюда и из непрерывности $\tl r$ имеем:
\[
  \tl r(2\,\pi) = \lim\limits_{t \to 2\pi-0} \tl r(t) = \lim\limits_{t \to 2\pi-0} x_5(t) = {}-2\, \pi.
\]

Теперь рассмотрим $\tl r([2\,\pi, 4\,\pi))$. Из определений $\tl r$ и $x(j,\cdot)$ вытекает, что
\begin{gather*}
  \tl r([2\,\pi, 4\,\pi)) \subseteq X_2 := x(\NN_5 \ts [2\,\pi, 4\,\pi)) = \bigcup_{k=1}^5 x_k[2\,\pi, 4\,\pi) = \{(0,+2)\} \cup \{(0,0)\} \cup \{(0,-2)\} \cup{}
\\
  {}\cup \{(x,y) \in \RR^2 \mid (0 \le x < 2\,\pi) \wedge (y = 0)\} \cup \{(x,y) \in \RR^2 \mid (-2\, \pi \le x < 0) \wedge (y = 0)\}.
\end{gather*}
При непрерывном отображении образ линейно связного множества (в данном случае полуотрезка $[2\,\pi, 4\,\pi)$) линейно связен, поэтому $\tl r([2\,\pi, 4\,\pi))$ лежит в одной из четырёх компонент линейной связности множества $X_2$. Но поскольку $\tl r(2\,\pi) = -2\,\pi$, то этой компонентой будет $x_5[2\,\pi, 4\,\pi)$, содержащее точку $(-2\,\pi,0)$. Ввиду того, что $x_5([2\,\pi, 4\,\pi))$ не пересекается с другими $x_k([2\,\pi, 4\,\pi))$, это значит, что $\forall t \in [2\,\pi, 4\,\pi)$ $\tl r(t) = x_5(t)$. Тогда из непрерывности $\tl r$ имеем:
\[
  \tl r(4\,\pi) = \lim\limits_{t \to 4\pi-0} \tl r(t) = \lim\limits_{t \to 4\pi-0} x_5(t) = 0.
\]

Рассматривая $\tl r((4\,\pi, 6\,\pi))$ и рассуждая также как и при рассмотрении $\tl r((0, 2\,\pi))$, приходим к выводу, что $\forall t \in (4\,\pi, 6\,\pi)$ $\tl r(t) = x_4(t)$. Отсюда и из непрерывности $\tl r$ имеем:
\[
  \tl r(4\,\pi) = \lim\limits_{t \to 4\pi+0} \tl r(t) = \lim\limits_{t \to 4\pi+0} x_4(t) = {}+2\, \pi.
\]
Но выше было установлено, что $\tl r(4\,\pi) = 0$. Мы пришли к противоречию. $\not\!\! \triangle$

\bibliography{Polynomials}

\end{document}